\documentclass[12pt]{amsart}

\usepackage{amsmath, amssymb, latexsym,amscd}
\begin{document}
\newtheorem{theorem}{Theorem}[section]
\newtheorem{lemma}[theorem]{Lemma}
\newtheorem{proposition}[theorem]{Proposition}
\newtheorem{fact}[theorem]{Fact}
\newtheorem{condition}[theorem]{Condition}
\newtheorem{example}[theorem]{Example}
\newtheorem{claim}[theorem]{Claim}
\newtheorem{question}[theorem]{Question}
\newtheorem{corollary}[theorem]{Corollary} 
\theoremstyle{definition}
\newtheorem{definition}[theorem]{Definition}
\newtheorem{statement}[theorem]{Statement}
\newtheorem{notation}[theorem]{Notation} 
\theoremstyle{remark}
\newtheorem{remark}[theorem]{Remark}
\newcommand\cl{\begin{claim}}
\newcommand\ecl{\end{claim}}
\newcommand\rem{\begin{remark}\upshape}
\newcommand\erem{\end{remark}}
\newcommand\ex{\begin{example}\upshape}
\newcommand\eex{\end{example}}
\newcommand\nota{\begin{notation}\upshape}
\newcommand\enota{\end{notation}}
\newcommand\dfn{\begin{definition}\upshape}
\newcommand\edfn{\end{definition}}
\newcommand\cor{\begin{corollary}}
\newcommand\ecor{\end{corollary}}
\newcommand\thm{\begin{theorem}}
\newcommand\ethm{\end{theorem}}    
\newcommand\prop{\begin{proposition}}
\newcommand\eprop{\end{proposition}}
\newcommand\lem{\begin{lemma}}
\newcommand\elem{\end{lemma}}
\newcommand\fct{\begin{fact}}
\newcommand\efct{\end{fact}}

\providecommand\qed{\hfill$\quad\Box$}
\newcommand\pr{{Proof:\;}}
\newcommand\prc{\par\noindent{\em Proof of Claim: }}
\newcommand\dom{{\text{dom}}}
\newcommand\Pn{{P_n}}
\newcommand\ord{{ord}}
\newcommand\rk{{rk}}
\newcommand\car{{\rm{char}}}
\newcommand\M{{{\mathcal M}}}
\newcommand\N{{{\mathcal N}}}
\newcommand\K{{{\mathcal K}}}
\newcommand\D{{{\mathcal D}}}
\newcommand\A{{{\mathcal A}}}
\newcommand\Ps{{{\mathcal P}}}
\newcommand\E{{\mathcal E}}
\newcommand\Y{{{\mathbf{Y}}}}
\renewcommand\O{{{\mathcal O}}}
\renewcommand\L{{{\mathcal L}}}
\newcommand\Th{{\text{Th}}}
\newcommand\IZ{{\mathbb Z}}
\newcommand\IQ{{\mathbb Q}}
\newcommand\IR{{\mathbb R}}
\newcommand\IN{{\mathbb N}}
\newcommand\IC{{\mathbb C}}
\newcommand\F{{\mathcal F}}
\newcommand\IF{{\mathbb F}}
\newcommand\Se{{\mathcal S}}
\newcommand\V{{\mathcal V}}
\newcommand\W{{\mathcal W}}
\newcommand\T{{\mathcal T}}
\newcommand\si{{\sigma}}
\newcommand\n{{\nabla}}
\newcommand\C{{\mathcal C}}
\newcommand\Lr{{{\mathcal L}_{\text{rings}}}}
\newcommand\Lrd{{{\mathcal L}^*_{\text{rings}}}}
\newcommand\B{{\mathcal B}}
\newcommand\La{{\mathcal L}}
\def\U{{ \mathfrak U}}
\def\B{{\mathcal B}}
\def\I{{\mathcal I}}
\def\R{{\mathbb R}}
\def\G{{\mathfrak G}}
\let\le=\leqslant
\let\ge=\geqslant
\let\subset=\subseteq
\let\supset=\supseteq

\newcommand{\bigdoublewedge}{%
  \mathop{
    \mathchoice{\bigwedge\mkern-15mu\bigwedge}
               {\bigwedge\mkern-12.5mu\bigwedge}
               {\bigwedge\mkern-12.5mu\bigwedge}
               {\bigwedge\mkern-11mu\bigwedge}
    }
}

\newcommand{\bigdoublevee}{%
  \mathop{
    \mathchoice{\bigvee\mkern-15mu\bigvee}
               {\bigvee\mkern-12.5mu\bigvee}
               {\bigvee\mkern-12.5mu\bigvee}
               {\bigvee\mkern-11mu\bigvee}
    }
}

 \title{On differential Galois groups of strongly normal extensions}
\date{\today}

 \author{Quentin Brouette}
\email{Quentin.Brouette@umons.ac.be}
  
  \address{D\'epartement de Math\'ematique (Le Pentagone)\\
Universit\'e de Mons\\
20 place du Parc\\ B-7000 Mons\\ Belgium}

  \author{Fran\c coise Point}
  \email{point@math.univ-paris-diderot.fr}
  \address{Research Director at the FRS-FNRS\\D\'epartement de Math\'ematique (Le Pentagone)\\
Universit\'e de Mons\\
20 place du Parc\\ B-7000 Mons\\ Belgium}

\begin{abstract}
We revisit Kolchin's results on definability of differential Galois groups of strongly normal extensions, in the case where the field of constants is not necessarily algebraically closed. 
In certain classes of differential topological fields, which encompasses ordered or p-valued differential fields, we find a partial Galois correspondence and we show one cannot expect more in general. 
In the class of ordered differential fields, using elimination of imaginaries in $CODF$, we establish a relative Galois correspondence for relatively definable subgroups of the group of differential order automorphisms.
\end{abstract}
\maketitle

\section{Introduction}
\par  We investigate Galois theory for strongly normal extensions within given classes
of differential fields of characteristic $0$, whose fields of constants
are not necessarily algebraically closed.
\par Our main motivation comes from inductive classes $\C$ of topological $\L$-fields
$K$ of characteristic $0$, as introduced in \cite{GP}, and their
expansion to differential structures where the derivation $D$ on $K$
has a priori no interactions with the topology on $K$ (see section 4.1). Let $\C^{ec}$ be the subclass of existentially closed elements of $\C$.
On $\C$, we assume the topology is first-order definable,  the class 
$\C^{ec}$ is elementary and Hypothesis $(I)$ holds. This hypothesis was introduced in \cite{GP} and says that any element of $\C$ has always an extension belonging to $\C$ which contains the field of Laurent series. In particular this implies that the elements of $\C^{ec}$ are large fields. 
Then, one can show that the corresponding class $\C_{D}$ of differential expansions of elements of $\C$ has a model-completion which is axiomatised by expressing that any differential polynomial in one variable with non zero separant, has a zero close to a zero of its associated algebraic polynomial \cite[Definition 5(3)]{GP}.   
\par The setting described above is an adaptation to a topological setting of a former work of M. Tressl \cite{T}. Recently N. Solanki described an alternative approach \cite{Sol} in a similar topological context.

\par Examples of such classes $\C$ are the classes of ordered fields, ordered valued fields, $p$-adic valued fields and valued fields of characteristic zero, and so the corresponding classes $\C^{ec}$ are the classes of real closed fields, real closed valued fields, $p$-adically closed fields and algebraically closed valued fields. When $\C$ is the class of ordered fields, we obtain in this way an axiomatisation of the model completion of $\C_{D}$, the class of closed ordered differential fields, 
alternative to the one given by 
M. Singer \cite{S}; the theory of closed ordered differential fields is denoted by $CODF$. 
\par Given a differential field $F$, we denote by $C_{F}$ its subfield of constants. Let $K\subseteq L$ be two differential fields. Let $\A$ be a saturated differentially closed field containing $L$ and let $\U$ be any intermediate differential field between $L$ and $\A$. 
\par Let $L/K$ be a strongly normal extension as defined by E. Kolchin \cite{K73}, where $C_K$ is not supposed to be algebraically closed. The field $L$ is first assumed to be finitely generated but in the last section we will drop that hypothesis.
\par We consider the differential Galois group $Gal(L/K)$, namely the group of differential field automorphisms of  $\langle L, C_{\A}\rangle$ fixing $\langle K, C_{\A}\rangle$, as well as, varying $\U$, the group $Gal_{\U}(L/K)$ of differential field automorphisms of $\langle L, C_{\U}\rangle$ fixing $\langle K, C_{\U}\rangle$; when $\U=L$ (and $C_L=C_K$), we denote it by $gal(L/K)$.
Kolchin showed that $Gal(L/K)$ is isomorphic to $H(C_{\A})$, for $H$ an algebraic group defined over $C_K$. Furthermore, it follows from his work \cite[chapter 6, section 4]{K73} that $Gal_{\U}(L/K)$ is isomorphic to $H(C_{\U})$.
We give a model-theoretic proof of that fact. Independently of our work, a more concise model-theoretic presentation is proposed in \cite{PSa}. 
\par We obtain the following partial Galois correspondence. 
Let $E$ be an intermediate extension $(K\subset E\subset L)$ such that $dcl^{\U}(E)\cap L=E$. 
We get that the map $E\mapsto Gal_{\U}(L/E)$ is injective but we observe that not every definable subgroup of $Gal_{\U}(L/K)$ defined over $C_K$ is of the form $Gal_{\U}(L/E)$ for some intermediate field $E$. 
\par We illustrate our results by revisiting the classical examples of strongly normal extensions, namely Picard-Vessiot and Weierstrass extensions and considering in each case the groups $Gal_{\U}(L/K)$.  
\par The existence and uniqueness of Picard-Vessiot extensions within specific
classes of differential fields (such as formally real or formally
p-adic differential fields) have been considered by several authors
(see for instance \cite{GGO}, \cite{CHP}, using among other things Deligne's work on Tannakian categories). It has been generalised to
strongly normal extensions corresponding to logarithmic differential
equations in \cite{KP}; there, they replace Tannakian categories by the formalism of groupoids as introduced by E. Hrushovski \cite{H}. 
\par In the case of a Picard-Vessiot extension, it is well-known \cite{Kap} that its differential Galois group is isomorphic to a linear algebraic  group defined over the field of constants. 
\par In the case of a Weierstrass extension of a formally real or a
formally $p$-adic field, we give more explicit information on
which groups may arise as $Gal_{\U}(L/K)$. When the field of constants is $\IR$ and the differential Galois group $G$ is infinite, we show that it is isomorphic either to $S^{1}$ or $S^{1}\times \IZ/2\IZ$, where $S^{1}$ is the compact connected real Lie group of dimension $1$. When the field of constants is an $\aleph_{1}$-saturated real closed field, we quotient out by $G^{00}$, the smallest type-definable subgroup of $G$ of bounded index \cite{HPP}, and we similarly obtain that $G/G^{00}$ is isomorphic to $S^1$ or to $S^{1}\times \IZ/2\IZ$. When the field of constants is $\IQ_{p}$, we show that $G$ is isomorphic to a subgroup of the $\IQ_{p}$-points of the associated elliptic curve.
\par Assuming now that $T$ is an axiomatisation of a class $\C$ of topological $\L$-fields where Hypothesis $(I)$ holds and letting $T_{c,D}^*$ an axiomatisation of the model-completion of $\C_{D}$, we easily obtain a criterion for the existence of a strongly normal extension within the class of models of $T_{D}$. We compare that result with Theorem 1.5 of \cite{KP}.
\par Then we restrict ourselves to the class of ordered differential fields. 
Let $L/K$ be two differential ordered fields, let $L\subset \U$ with $\U\models CODF$; let $Aut_{T}(L/K)$ be the subgroup of $Gal_{\U}(L/K)$ of automorphisms that preserve the order. Note that it is not a priori a definable subgroup (but simply an infinitely definable one with parameters in $L$). Let $G$ be a subgroup of $Gal_{\U}(L/K)$, let $\mu(G)$ be its image in $H(C_{\U})$ and assume that $\mu(G)$ is definable with parameters in $\langle L,C_{\U}\rangle$. 
Then, using that the theory $CODF$ has elimination of imaginaries \cite{Point}, we show that 
we can find a tuple $\bar d$ in the real closure of $\langle L,C_{\U}\rangle$ such that $Aut_{T}(L\langle \bar d\rangle/K\langle \bar d\rangle)$ is isomorphic to $G\cap Aut_{T}(L/K)$. 
\par In the last section we consider non necessarily finitely generated strongly normal extensions $L/K$, as defined by J. Kovacic \cite{Kov1}, and show
that $Gal_{\U}(L/K)$ is isomorphic to a subgroup of an inverse limit of groups of the form $H_i(C_{\U})$, where $H_i$ is an algebraic group defined over $C_{K}$, $i\in I$.
\par This article supersedes the results obtained in \cite{BP} where we dealt with the case of ordered differential fields and are also a part of the first author's thesis \cite[chapter 3]{B}. 
 \section{Preliminaries}\label{basics}


\par Throughout the paper, $\L$ denotes a first order language, $\M$ an $\L$-structure and $A\subset M$. Only in this section, we will make the distinction between an $\L$-structure $\M$ and its domain $M$.
\nota We denote by $\L(A):=\L\cup\{c_{a}: a\in A\}$, where $c_{a}$ is a new constant symbol for each element $a\in A$ 
(not to be confused with the element of the subfield of constants in a differential field).  
 
\par Given a tuple $\bar a\in M$, we denote by $tp^{\M}(\bar a/A)$
the {\it type} of $\bar a$ in $\M$ over $A$.
\par  The set of all algebraic elements over $A$ in $\M$ is denoted by $acl^{\M}(A)$ and called the (model-theoretic) algebraic closure of $A$ in $\M$.
Similarly the set of definable elements over $A$ in $\M$ is denoted by $dcl^{\M}(A)$
 and called the definable closure of $A$ in $\M$.
 \enota
\par We will consider (expansions of) differential fields as $\La_{D}$-structures where $\L$ contains the language of rings $\L_{rings}:=\{+,-,.,0,1\}$ and $\La_{D}:=\La \cup\{^{-1}, D\}$ where $D$ is a new unary function symbol that will be interpreted in our structures as a derivation:
\begin{description}
	\item[Additivity] $\forall a\,\forall b \;\;D(a+b)=D(a)+D(b),$
	\item[Leibnitz rule] $\forall a\,\forall b\;\;\; D(a.b)=a.D(b)+D(a).b .$
\end{description}
\par We denote by $DCF_{0}$ (respectively $ACF_{0}$) the theory of differentially (respectively algebraically) closed fields of characteristic $0$. 
\par A. Robinson showed that $DCF_{0}$ (respectively $ACF_0$) is the model completion of the $\La_{rings,D} \cup\{^{-1}\}$-theory of differential fields (of, respectively, the $\La_{rings} \cup\{^{-1}\}$-theory of fields) of characteristic $0$. Since these theories can be universally axiomatised, it implies that $DCF_{0}$ (respectively $ACF_0$) admits quantifier elimination \cite[Theorem 13.2]{Sa}. 
\par We will use that $ACF_{0}$ and $DCF_{0}$ admit elimination of imaginaries (e.i.) (the latter is deduced from the former) \cite[Section 8.4]{TZ}. This property of the theory $ACF_{0}$ can be seen as the model-theoretic counterpart of the fact that any ideal has a field of definition;
it can also be seen more concretely as follows.
 Let $F_{0}\subset F$, with $F\models ACF_{0}$ and $F_{0}$ a subfield, recall that given an $F_{0}$-definable function $f$ in $F$, we can find finitely many $F_{0}$-algebraic subsets $X_{i}$ such that $f$ is equal to an $F_{0}$-rational function on $X_{i}$ (straightforward from \cite[Proposition 3.2.14]{Ma}).  
Moreover if $E$ is an $F_0$-definable equivalence relation on $F^{n}$, then there is an $F_{0}$-definable function $f$ from $F^{n}$ to $F^{\ell}$ 
such that $E(\bar x,\bar y)$ iff $f(\bar x)=f(\bar y)$ \cite[Theorem 3.2.20]{Ma}. 
\smallskip
\par L. Blum gave an elegant axiomatisation of $DCF_0$ \cite[Section 40]{Sa} and showed that any differential field $K$
has an atomic extension, model of $DCF_0$, called the differential closure of $K$ \cite[Section 41]{Sa}. Recall that it is unique up to isomorphism and that the type of a tuple $\bar a$ in the differential closure of $K$ is isolated over $K$; 
following E. Kolchin, 
one says that $\bar a$ is {\it constrained} over $K$ \cite[Section 2]{K74}.  
In the introduction of \cite[page 142]{K74}, 
one can find a "dictionary" between the terminology used by model theorists and by the differential algebra community. 
\medskip



\nota \label{prem} 
Given a field $F$ (respectively a differential field $F$), denote by $\bar F$ the algebraic closure (respectively $\hat{F}$ some differential closure) of $F$. 
We denote by $C_{F}$ the subfield of constants of $F$.
Note that $C_{\hat{F}}=\bar C_{F}$.
\par Given an algebraic group $H$ defined over $F_0\subset F$ and given a subfield $F_0\subset F_1\subset F$, we will denote by $H(F_1)$ its $F_1$-points.
\enota
\par Now let us introduce some differential algebra standard notation.
\nota\label{dfninitial} Let $K\{X_{1},\cdots,X_{n}\}$ be the  ring of differential polynomials over $K$ in $n$ differential indeterminates $X_{1},\cdots, X_{n}$ over $K$, namely it is the ordinary polynomial ring in indeterminates $X_{i}^{(j)}$, $1\leq i\leq n$, $j\in \omega$, with by convention $X_{i}^{(0)}:=X_{i}$. One can extend the derivation $D$ of $K$ to this ring by setting $D(X_{i}^{(j)}):=X_{i}^{(j+1)}$ and using the additivity and the Leibnitz rule. 

\par Set $\bar X:=(X_{1},\cdots, X_{n})$ and $\bar X^{(j)}:=(X_{1}^{(j)},\cdots, X_{n}^{(j)})$. 
We will denote by $K\langle \bar X\rangle$ the fraction field of $K\{\bar X\}$.
Let $f(\bar X)\in K\{\bar X\}\setminus K$ and suppose that $f$ is of order $m$, then we associate with 
$f(\bar X)$ the ordinary polynomial $f^*(\bar X^*)\in K[\bar X^*]$ with $\bar X^*$ a tuple of variables of the appropriate length such that $f(\bar X)=f^*(\bar X^{\n})$ where $\bar X^{\n}:=(\bar X^{(0)},\ldots ,\bar X^{(m)})$.
 We will make the following abuse of notation: if $\bar b\in K^n$, then $f^*(\bar b)$ means that we evaluate the polynomial $f^*$ at the tuple $\bar b^{\n}$.

\par If $n=1$, recall that the separant $s_{f}$ of $f$ is defined as $s_{f}:=\frac{\partial f }{\partial X_{1}^{(m)}}$.
\enota
\medskip
\par We will use the following consequence of the quantifier elimination results for $T=ACF_{0}$ (respectively $DCF_{0}$). 
Let $\M\models T$, let $K$ be a subfield of $M$ (respectively a differential subfield) and $\bar a\in M$, 
then the type $tp^{\M}(\bar a/K)$ of $\bar a$ over $K$ is determined by the ideal $$\I_{K}(\bar a):=\{p(\bar X)\in K[\bar X]:\;p(\bar a)=0\}$$ (respectively the differential ideal $\I^D_{K}(\bar a):=\{p(\bar X)\in K\{\bar X\}:\;p(\bar a)=0\}$).
\medskip
\par We will work inside a {\it saturated} structure $\M$ \cite[Definition 4.3.1]{Ma}, but as usual we will only need $\kappa^+$-saturation for some cardinal $\kappa$.
Such a structure $\M$ has the following {\it homogeneity} property \cite[Proposition 4.3.3 and Proposition 4.2.13]{Ma}.
Let $A\subset M$ such that $|A|<|M|$ and $\bar a, \bar b$ be finite tuples of elements of $M$ 
with $tp^{\M}(\bar a/A)=tp^{\M}(\bar b/A)$, 
then there is an automorphism of $\M$ fixing $A$ and sending $\bar a$ to $\bar b$.
\medskip

\section{Strongly normal extensions and differential Galois groups}

\par This is mainly an expository section on differential Galois groups.
\par Let $K\subseteq L\subseteq \U$ be three differential fields and let $C_K\subset C_L\subset C_{\U}$ their respective subfields of constants; we will also use the notations $\U/L/K$, $L/K$ or $\U/L$.  Let $\A$ be a saturated model of $DCF_0$ containing $\U$ of cardinality strictly bigger than $\vert K\vert$; $\A$ will play the role of a universal extension \cite[page 768]{K53}. 
\par Since $DCF_{0}$ admits quantifier elimination, $\A$ has the following property: any $\L_{rings,D}$-embedding from $K\langle \bar a \rangle$ to $K\langle \bar b \rangle $ fixing $K$ and sending $\bar a$ to $\bar b$, where $\bar a, \bar b$ are finite tuples from $\A$, can be extended to an automorphism of $\A$ fixing $K$ (we use the homogeneity of $\A$, see section \ref{basics}). 
\par Denote by $Hom_{K}(L,\A)$ the set of $\L_{rings,D}$-embeddings from $L$ to $\A$ which are fixed on $K$. 
\medskip

\par The notion of {\it strongly normal} extension has been introduced by E. Kolchin \cite[page 393]{K73}.  
\dfn An element $\tau\in Hom_{K}(L,\A)$ is {\it strong} if $\tau$ is the identity on $C_{L}$ and if $\langle L,C_{\A} \rangle=\langle\tau(L),C_{\A} \rangle$  \cite[page 389]{K73}.
The extension $L/K$ is a {\it strongly normal extension} if $L$ a finitely generated extension of $K$ such that every $\tau\in Hom_{K}(L,\A)$ is {\it strong}. 
\edfn
\par From the above remarks, it follows that, whenever $L$ is a finitely generated extension of $K$, any element of $Hom_{K}(L,\A)$ extends to an $\L_{rings,D}$-automorphism of $\A$ fixing $K$ and when $L$ is a strongly normal extension of $K$, any such automorphism restricts to an automorphism of $\langle L,C_{\A} \rangle.$
\par One can extend the notion of ''strongly normal'' to not necessarily finitely generated extensions $L$ of $K$ by asking that $L$ is the union of finitely generated strongly normal extensions \cite[Chapter 2, section 1]{Kov1}. This will be done in Section \ref{nfg}.
\fct \cite[Chapter 6, section 3, Proposition 9]{K73} If $L$ is a strongly normal extension of $K$, then $C_{K}=C_{L}$. 
\efct
\nota \label{auto-nota} Let $\U/L/K$ be two differential extensions of $K$. We denote by $gal(L/K)$ the group of $\L_{rings,D}$-automorphisms of $L$ fixing $K$, by $Gal(L/K)$ the group of $\L_{rings,D}$-automorphisms of $\langle L, C_{\A}\rangle$ fixing $\langle K, C_{\A}\rangle$, and by $Gal_{\U}(L/K)$ the group of $\L_{rings,D}$-automorphisms of $\langle L, C_{\U}\rangle$ fixing $\langle K, C_{\U}\rangle$.
\enota
\medskip
\par Let $L$ be a strongly normal extension of $K$.
We will revisit Kolchin's result that $Gal(L/K)$ is isomorphic to the $C_{\A}$-points of an algebraic group $H$ defined over $C_{K}$ and in particular we show that $Gal_{\U}(L/K)$ is isomorphic to the $C_{\U}$-points of $H$ for any field $L\subset \U\subset \A$ (Theorem \ref{Kolchin}). 


\par We begin by a few easy remarks, using that the theory $DCF_0$ admits quantifier elimination and basic facts about linear disjointness \cite[Chapter 3, Theorem 2]{Lang}. 
\rem\label{ld} Suppose that $C_{K}$ is algebraically closed in $C_{\U}$ (denoted by $C_K\subset_{ac} C_{\U}$). 
Then the fields $\U$ and $\bar C_{K}$ are linearly disjoint over $C_{K}$, respectively $\U$ and $\langle K,\bar C_{K}\rangle$ are linearly disjoint over $K$. Therefore $\langle L,\bar C_{K}\rangle \cap \U=L$. 
\par Moreover, $\U$ (respectively $L$) and $C_{\A}$ are linearly disjoint over $C_{\U}$ (respectively $C_{K}$) and so $\langle L, C_{\A}\rangle \cap \U=\langle L,C_{\U}\rangle$.\qed
\erem
\rem\label{star} Assume now that $L/K$ is a strongly normal extension. 
By Remark \ref{ld}, for any $\si\in Hom_{K}(L,\U)$, we have $\langle L, C_{\U}\rangle =\langle \si(L), C_{\U}\rangle$.
Moreover, letting $\bar u\in L$, we get: $tp^{\A}(\bar u/K)\models tp^{\A}(\bar u/\langle K,\bar C_{K}\rangle)$.
 \par This follows from the fact that $C_L=C_K$, and so $L$ and $\langle K,\bar C_{K}\rangle$ are linearly disjoint over $K$, which implies that the ideal $\I^D_{\langle K,\bar C_{K}\rangle}(\bar u)$ (Notation \ref{prem}) has a set of generators in $K\{\bar X\}$ \cite[Chapter 3, Section 2, Theorem 8]{Lang}.
 \qed
\erem
\medskip
\fct\label{Tsn}  \cite[section 9, Theorem 3]{K74}, \cite{P98} Let $L=K\langle \bar a \rangle $ be a strongly normal extension of $K$. Then there is a formula $\xi(\bar y)\in tp(\bar a/\langle K, C_{\A} \rangle)$ such that for any $\bar b\in \A$ with $\A\models \xi(\bar b)$, we have $\bar a\in dcl(K, \bar b,C_{\A})$. Consequently by choosing $\bar b$ in a differential closure of $\langle K, C_{\A} \rangle$, we get that $tp(\bar a/\langle K, C_{\A} \rangle)$ is isolated.
\efct
\par Therefore, since $\hat{K}$ embeds in any model of $DCF_0$ extending $K$, the definition of {\it strongly normal} does not depend on the universal extension one is working in.  Moreover one can observe that ''$K\langle \bar a\rangle$ is a strongly normal extension of $K$'' only depends on $tp^{\A}(\bar a/K)$, see \cite{P98}. 
\medskip
\par From now on, for ease of notation, we will denote the strongly normal extension $L$ by $K\langle a\rangle$, not making the distinction whether $a$ is a single element or a tuple of elements (as it is customary done in model theory). We will extend this abuse of notation to tuples $\bar f$ of definable functions applied to such tuple $\bar a$, by denoting it by $f$. It is straightforward to check it has no incidence on the proofs. Besides, we could also have used a recent result of Pogudin \cite{Pog} (however we only use it in the proof of Lemma \ref{KP}).

\rem\label{ei} Let $L:=K\langle  a\rangle$. By the above fact, there is an $\L_{rings,D}(K)$-formula $\phi( y,\bar c)$ with the tuple $\bar c\in C_{\A}$, isolating $tp^{\A}( a/\langle K,C_{\A}\rangle)$. Since $DCF_{0}$ admits elimination of imaginaries, we may choose $\bar c$ to be a canonical parameter for the definable set $\phi(\A,\bar c)$ \cite[Definition 8.4.1]{TZ} and so we may assume that $\bar c\in dcl^{\A}(K\langle  a\rangle)$ (and so $\bar c\in L\cap C_{\A}=C_K$). By abuse of notation, we keep the same formula $\phi$ and the tuple $\bar c$. Since now $\bar c\in C_{K}$ and $\phi( y,\bar c)$ isolates $tp^{\A}( a/\langle K,C_{\A}\rangle)$, for any intermediate field $L\subset \U\subset \A$, it also isolates $tp^{\A}( a/\langle K,C_{\U}\rangle)$.
\erem
\par Note that for $ a\in L$, in order to show that $tp^{\A}( a/K)$ implies $tp^{\A}( a/\langle K,C_{\A}\rangle)$, one can also use that $L$ and $C_{\A}$ are linearly disjoint over $C_{K}$. 

\lem\label{emb}  Let $L/K$ be a strongly normal extension, then there is an embedding of $gal(L/K)$ into $Gal(L/K)$. 
Any element of $Hom_{K} (L,\U)$ extends uniquely to an element of $Gal(L/K)$, which restricts to an element of $Gal_{\U}(L/K)$.
\elem
\pr Let $ a\in L$ be such that $L=K\langle  a \rangle$.
\par Let $\si\in Hom_{K} (L,\U)$. Since $DCF_{0}$ admits q.e., $tp^{\A}(a/K)=tp^{\A}(\si( a)/K)$.

So, using Remark \ref{ei}, we get that $tp^{\A}(a/\langle K, C_{\A}\rangle)=tp^{\A}(\si(a)/\langle K, C_{\A}\rangle).$ 
Thus the map $\tilde \si$ sending $a$ to $\si(a)$ and which is the identity on $\langle K, C_{\A}\rangle$ is elementary. Going to a $\vert C_{\A}\vert^+$-saturated extension $\tilde \A$ of $\A$, we can extend it to an automorphism of that extension and 
since $L$ is strongly normal, 
we have $\langle L,C_{\tilde \A}\rangle=\langle K,\si(a), C_{\tilde \A}\rangle$. By Remark \ref{ld}, $\langle L,C_{\tilde \A}\rangle\cap \U=\langle L,C_{\U}\rangle$ and $\langle \si(L),C_{\tilde \A}\rangle\cap \U=\langle \si(L),C_{\U}\rangle$. So, $\tilde \si\in Gal_{\U}(L/K)$. The uniqueness is clear. 
\par In the case when $\U=L$, we get that $gal(L/K)$ embeds in $Gal(L/K)$. 
\qed
\medskip

\par Now let us state the following result which follows from former work of Kolchin {\rm \cite{K73}}.
\thm \label{Kolchin}  Let $L/K$ be a strongly normal extension. Then, the group $Gal(L/K)$ is isomorphic to the $C_{\A}$-points of an algebraic group $H$ defined over $C_K$.
Moreover for any intermediate subfield $L\subset \U\subset \A$, $Gal_{\U}(L/K)$ is isomorphic to $H(C_{\U})$. In particular, $gal(L/K)$ is isomorphic to $H(C_{K})$.
\ethm 
\pr Let $\phi(y,\bar c)$, $\bar c\in C_{K}$, be the formula obtained in Remark \ref{ei}, isolating $tp^{\A}( a/\langle K,C_{\A}\rangle)$.
By the compactness theorem (using the fact that for any two realisations $a,\;b$ of $\phi(y,\bar c)$, $a\in dcl(K,\bar b,C_{\A})$), there is an $\L_{rings,D}(K)$-definable function $f(.,.)$ such that 
given $b_{1},\;b_{2}\in \phi(\A,\bar c)$, there is $\bar c_{1,2}\in C_{\A}$ such that $f(b_{1},\bar c_{1,2})=b_{2}$ $(\star)$. Again applying the elimination of imaginaries for $DCF_{0}$, we may assume in $(\star)$ that $\bar c_{1,2}$ is unique. Therefore given $\si\in Gal(L/K)$, we can associate a unique tuple $\bar c_{a,\si(a)}$ in $C_{\A}$ such that $f( a,\bar c_{a,\si(a)})=\si(a)$. So the map from $Gal(L/K)$ to some cartesian product $C_{\A}^n$ sending $\si$ to $\bar c_{a,\si(a)}$ is injective, where $n$ is equal to the length of $\bar c_{a,\si(a)}$. 
Re-write $f(.,.)$ as $\tilde f(.,.,\bar k)$ with $\tilde f$ a $\L_{rings,D}$-definable function.
\par Let $S:=\{\bar d\in C_{\A}:\;\exists b_{1}\exists b_{2}\;\bigwedge_{i=1}^2   \phi(b_{i},\bar c)\;\&\;\tilde f(b_{1},\bar d,\bar k)=b_{2}\}$.
Now we use that the type (in $DCF_{0}$) of $\bar k$ over $C_{\A}$ is definable, to get an $\L_{rings}(C_{\A})$-definition of $S$. In other words, $S$ is stably embedded. Since $S$ is $K$-definable, we may conclude that $S$ is in fact $C_K$-definable. (A detailed argument of that fact is the following. Use elimination of imaginaries in $ACF_{0}$ and the fact that $S$ is invariant by all automorphisms of $\A$ which are fixed on $K$, in order to deduce that there is a canonical parameter for $S$ belonging to $C_{\A}\cap K=C_K$.)
\par One can define a group law $\otimes$ on $S$ as follows: let $\bar c_1,\;\bar c_2\in S$, then $\si_1(a)=\bar f(a,\bar c_1)$ and $\si_2(a)=f(a,\bar c_2)$. So there is a unique tuple $\bar c_3$ associated with $\si_3:=\si_2\circ \si_1$, namely $f(a,\bar c_3)=f(f(a,\bar c_1),\bar c_2)$ and we define $\bar c_3:=\bar c_1\otimes \bar c_2$. It is defined by the formula: $\exists b\;(\phi(b,\bar c)\;\&\; f(b,\bar c_3)=f(f(b,\bar c_1),\bar c_2))$. Similarly as above, one shows that the group law is $\L_{rings}(C_{K})$-definable.
\par Then we use the fact that a $\L_{rings}(C_{K})$-definable group in an algebraically closed field is an algebraic group defined over $C_{K}$. Therefore, $S=H(C_{\A})$, for some algebraic group $H$ defined over $C_K$.
\par Finally let us check that $H(C_{\U})$ is isomorphic to $Gal_{\U}(L/K)$. Let $\bar d\in C_{\U}$, then $\si(a)=f(a,\bar d)\in \langle L,C_{\U}\rangle$. Thus $\si(L)\subset \langle L,C_{\U}\rangle$ and so $\si\in Gal_{\U}(L/K)$. Conversely, let $\si\in Gal_{\U}(L/K)$, so $\si(a)\in \U$. Let $\bar d'\in C_{\A}$ be such that $f(a,\bar d')=\si(a)$ and since we have assumed that $\bar d'$ was uniquely determined, it belongs to $dcl(\langle L,C_{\U}\rangle)\cap C_{\A}=C_{\U}.$ 
\qed
\medskip
\nota \label{mu} We will denote by 
$\mu$ 
the isomorphism either from $Gal_{\U}(L/K)$ to $H(C_{\U})$ or its restriction from $gal(L/K)$ to $H(C_{K})$.
\enota
\cor \label{GalE} Let $L/K$ be a strongly normal extension and let $E$ an intermediate differential subfield.
Then $Gal_{\U}(L/E)$ is isomorphic to an algebraic subgroup of $H(C_{\U})$ defined over $C_K$.
\par If $E$ is a strongly normal extension of $K$, then $Gal_{\U}(L/E)$ is a normal subgroup of $Gal_{\U}(L/K)$.
\ecor
\pr Note that $L$ is a strongly normal extension of $E$ and that $E$ is finitely generated over $K$, say $E:=K\langle \bar e\rangle$ with $\bar e\subset L$ ($L$ is also finitely generated as a field extension of $K$ and so one applies \cite[Chapter 3, section 2, Proposition 6]{Lang}). 
\par We proceed as before noting that to express that $\si\in Gal_{\U}(L/E)$, it suffices to say that $\si\in Gal_{\U}(L/K)$ and $\si(\bar e)=\bar e$. 
We will make the same abuse of notation as before, denoting $\bar e$ by simply $e$.
So we express $e$ in terms of $a$ and $K$, namely $e=h(a)$ for some $h$ an $\L_{rings,D}(K)$-definable function. Using the same notation as in Theorem \ref{Kolchin}, we get that: 
\begin{equation}
\si(a)=\tilde f(a,\bar c_{a,\si(a)}, \bar k), \end{equation} 
\par Now we express that $\si(e)=e$ by:
\begin{equation}  h(a)=h(\tilde f( a,\bar c_{a,\si(a)}, \bar k)). \end{equation}
We consider the $\L_{rings,D}(K)$-definable set $S_{E}:=\{\bar d\in C_{\A}:\;\forall y\;(\phi( y,\bar c)\rightarrow h(y)=h(\tilde f(y,\bar d, \bar k))\}$.
Again, we use the fact that $tp^{\A}(\bar k/C_{\A})$ is definable and we get a $\L_{rings}(C_{\A})$-formula $\theta_{E}(\bar z)$ such that for any $\bar d\in C_{\A}$, we have $\theta_{E}(\bar d)$ holds if and only if $\forall y\;(\phi(y,\bar c)\rightarrow h(y)=h(\tilde f(y,\bar d, \bar k))$ belongs to $tp^{\A}(\bar k/C_{\A})$. By elimination of imaginaries for $DCF_{0}$, we get a canonical parameter for that formula belonging to $dcl^{DCF_0}(K)\cap C_{\A}=C_K$. Restricting down to $C_{\U}$, we obtain an $\L_{rings}(C_K)$-definable subgroup of $H(C_{\U})$.
\par Let us prove the second assertion. Let $\tau\in Gal_{\U}(L/K)$, $\si\in Gal_{\U}(L/E)$ and $b\in E$. We want to show that $\tau(b)$ belongs to $\langle E, C_{\U}\rangle$. Since $E$ is a strongly normal extension of $K$ and $\tau\in Hom_{K}(E,\U)$, by Remark \ref{star}, $\langle E,C_{\U}\rangle=\langle \tau(E),C_{\U}\rangle$.
\qed

\section{Normality and Galois correspondence}
\par In this section, we apply the results of the preceding section in order to deduce a partial Galois correspondence. In the next section, we will show that one cannot expect more in general.
\par We will first examine under which hypothesis, a strongly normal extension of $K$ included in $\U$ is $Gal_{\U}(L/K)$-normal, we use a similar terminology as used by I. Kaplansky for a differential extension to be normal \cite[page 20]{Kap}, namely:
\dfn\label{normal} Let $L$ be a differential field extension of $K$ with $K\subset L\subset \A$. Let $\Sigma\subset Hom_{K}(L,\A)$. Then $L$ is $\Sigma$-normal if for any $x\in L\setminus K$, there is $\tau\in \Sigma$ such that $\tau(x)\neq x$.
\edfn
\par I. Kaplansky showed that if $K$ is a differential field of characteristic $0$ with algebraically closed field of constants, then any Picard-Vessiot extension $L$ of $K$ is $gal(L/K)$-normal \cite[Theorem 5.7]{Kap}. 
E. Kolchin showed that if $L$ is a strongly normal extension of $K$, then $L$ is $Hom_{K}(L,\A)$-normal \cite[Chapter 6, section 4, Theorem 3]{K73}.
\medskip

\prop \label{normalext} Assume that $\U$ is $\vert K\vert^+$-saturated. Suppose that $L/K$ is a strongly normal extension and that $dcl^{\U}(K)\cap L=K$. Then $L$ is a $Gal_{\U}(L/K)$-normal extension of $K$.
\eprop
\pr Let $L:=K\langle a\rangle$. By Lemma \ref{emb}, it suffices to prove that $L$ is $Hom_{K}(L,\U)$-normal.
\par Let $u\in L\setminus K$ and consider $tp^{\U}(u/K)$. Express $u$ as a rational function $\frac{p_{1}(a)}{p_{2}(a)}$, with $p_{1}(X),\;p_{2}(X)\in K\{X\}$ and $p_{2}(a)\neq 0$. By Fact \ref{Tsn} and Remark \ref{ei}, $tp^{\A}(a/K)$ is isolated by a quantifier-free formula $\chi(x)$ with parameters in $K$.
\par Suppose there is $d\neq u\in \U$ such that $tp^{\U}(u/K)=tp^{\U}(d/K)$. In particular for some $b\in \U$, we have $d=\frac{p_{1}(b)}{p_{2}(b)}\;\&\;\chi(b)$.
In particular, $tp^{\A}(b/K)=tp^{\A}(a/K)$. So, there is an element $\si\in Hom_{K}(L,\U)$ such that $\si(a)=b$.
\par So, either there is $u\neq d\in \U\setminus K$ such that $tp^{\U}(u/K)=tp^{\U}(d/K)$ and so we get that for some $\tilde\si\in Gal_{\U}(L/K)$, $\tilde\si(u)=d$.
\par Or $u$ is the only element of $\U$ realizing $tp^{\U}(u/K)$, then $u\in dcl^{\U}(K)$ and so by hypothesis $u\in K$. \qed
\medskip
 \par From the results recalled in the previous section, we deduce a partial Galois correspondence, between intermediate extensions and algebraic subgroups of $H(C_{\U})$.
\cor \label{interm} Let $K\subset L\subset \U$, suppose $L/K$ is strongly normal. Then,
\begin{enumerate}
\item the map sending $E$ to $Gal_{\U}(L/E)$ is an injective map from
  the set of intermediate extensions $E$ with $dcl^{\U}(E)\cap L=E$,
  to the set of algebraic subgroups of $H(C_{\U})$ defined over $C_K$;
 \item assume that $K\subset E\subset L$ is a strongly normal extension of $K$, then \\$Gal_{\U}(L/K)/Gal_{\U}(L/E)$ is isomorphic to an algebraic subgroup of $Gal_{\U}(E/K)$ defined over $C_K$.
 \end{enumerate} 
\ecor
\pr 
\par $(1)$ It follows from Proposition \ref{normalext} and Corollary \ref{GalE}.
\par $(2)$ Let $f$ be the map sending $\si\in Gal_{\U}(L/K)$ to $\si\restriction\langle E,C_{\U}\rangle$; it is well-defined and it goes from $Gal_{\U}(L/K)$ to $Gal_{\U}(E/K)$ by Remark \ref{star}. 
The kernel of $f$ is equal to $Gal_{\U}(L/E)$. Since both groups $Gal(L/K)$ and $Gal(L/E)$ are isomorphic to algebraic groups defined over $C_{K}$ and
since $ACF_{0}$ admits elimination of imaginaries, the quotient $Gal(L/K)/Gal(L/E)$ is isomorphic to a group $G(C_{\A})$, where $G$ is an algebraic group defined over $C_{K}$. 
Since the definable function which chooses a point in each  $C_{K}$-definable equivalence class (containing a $C_{\U}$-point) is also $C_{K}$-definable (see section 2), we obtain that $Gal_{\U}(L/K)/Gal_{\U}(L/E)$ is isomorphic to a $C_K$-definable subgroup of $G(C_{\U})$.
\par Since we have a Galois correspondence when dealing with the full Galois group $Gal$, we have that $Gal(L/K)/Gal(L/E)$ is isomorphic to $Gal(E/K)$. Therefore, $Gal_{\U}(E/K)$ is isomorphic to $G(C_{\U})$.\qed
\medskip
\par We will show that $Gal_{\U}(L/K)/Gal_{\U}(L/E)$ is, in general, isomorphic to a proper subgroup of $Gal_{\U}(E/K)$ (see Example \ref{proper}).

\section{Examples of strongly normal extensions}
\par We will provide in this section some examples of strongly normal extensions and of their Galois groups within some classes of differential fields. In particular, we review the two classical examples of Picard-Vessiot and Weierstrass extensions in classes of formally real fields or formally p-adic fields.  In this setting (and more generally), deep results on existence of strongly normal extensions have been obtained for instance in \cite{CHP}, \cite{GGO} and \cite{KP}. 
\par First we will recall the framework developed in \cite{GP}.
Then, we will review classical examples of strongly normal extensions as Picard-Vessiot extensions and Weierstrass extensions in these classes of topological differential $\L$-fields. 
\subsection{Model completion} Let $\La$ be a relational expansion of $\L_{rings}$ by $\{R_i; i\in
I\}\cup\{c_j;j\in J\}$ where the $c_j$'s are constants and the $R_i$'s 
are $n_{i}$-ary predicates, $n_{i}>0$. 
\dfn {\rm \cite{GP}}
Let $K$ be an $\La\cup \{^{-1}\}$-structure such that its restriction to $\La_{\text{rings}}\cup \{^{-1}\}$ is a field of characteristic $0$. Let $\tau$ be a  Hausdorff topology on $K$. Recall that $\langle K,\tau\rangle$ is a topological $\La$-field if the ring operations are continuous, the inverse function is continuous on $K\setminus\{0\}$ and every relation $R_{i}$ (respectively its complement $\lnot R_{i}$), with $i\in I$, is interpreted
in $K$ as the union of an open set $O_{R_{i}}$
(respectively $O_{\lnot R_i}$) and an algebraic subset $\{\bar x\in
K^{n_{i}}:\bigwedge_k r_{i,k}(\bar x)=0\}$ of $K^{n_i}$ (respectively
$\{\bar x\in K^{n_{i}}:\bigwedge_{l} s_{i,l}(\bar x)=0\}$ of $K^{n_{i}}$), where $r_{i,k},\;s_{i,l}\in K[X_1,\cdots,X_{n_i}]$.
\edfn
\par The topology $\tau$ is said to be first-order definable \cite{P87} if
there is a formula $\phi(x,\bar y)$ such that the set of subsets of the form $\phi(K,\bar a):=\{x\in K:\;K\models \phi(x,\bar a)\}$ with $\bar a\in K$ 
can be chosen as a basis $\V$ of neighbourhoods of $0$ in $K$. 
\par Examples of topological $\L$-fields are given in \cite[Section 2]{GP}. For instance, ordered fields, ordered valued fields, valued fields, $p$-valued fields, fields endowed with several distinct valuations or several distinct orders.
\par We now consider the class of models of a universal theory $T$ which has a model-completion $T_{c}$, and we assume that the models of  $T$ are in addition topological $\La$-fields, with a first-order definable topology.
Further, we work under the extra assumption that the class of models of $T_{c}$ satisfies {\it Hypothesis
$(I)$} \cite[Definition 2.21]{GP}. This hypothesis is the analog in our topological setting of the notion of {\it large fields} introduced by Pop \cite{Pop}.
\par We then consider the expansions of the models of $T$ to $\La_{D}$-structures and we denote by $T_D$ the $\La_D$-theory consisting of $T$ together the axioms expressing that $D$ is a derivation.
\par Hypothesis $(I)$ was used in \cite[Proposition 3.9]{GP} in order to show that any model of $T_D$ embeds in a model of $T_{c}$ which satisfies the scheme $(DL)$, namely that  if for 
each differential polynomial
 $f(X)\in K\{X\}\setminus\{0\}$, with $f(X)=f^*(X,X^{(1)},\ldots ,X^{(n)})$,
for every $W\in \V$, 

$\displaystyle (\exists \alpha_0,\ldots ,\alpha_n\;\;
f^*(\alpha_0,\ldots ,\alpha_n)=0 \wedge s_f^*(\alpha_0,\ldots ,\alpha_n)\ne 0 )
\Rightarrow$\\
$\displaystyle \Big(\exists z\big(f(z)=0\wedge s_f(z)\ne 0\wedge$
$\displaystyle \bigwedge_{i=0}^n (
z^{(i)}-\alpha_i\in W)\big)\Big)$.
\par Note that in \cite{GP}, the scheme $(DL)$ is not quite given as above, but in an equivalent form \cite[Proposition 3.14]{GP}.
\medskip
\par Note that since we assumed that the topology is first-order definable, the scheme
of axioms $(DL)$ can be expressed in a first-order way. Let $T_{c,D}^*$ be the $\La_{D}$-theory consisting of $T_{c}\cup T_{D}$ together with the scheme $(DL)$. Then $T_{c,D}^*$ is the model completion of $T_{D}$. A consequence of that axiomatisation is that the subfield of constants of a model $\U$ of $T_{c,D}^*$ is dense in $\U$ \cite[Corollary 3.13]{GP}.
\medskip
\nota Let $f\in K\{X\}$ and let $\langle f\rangle$ denote the differential ideal generated by $f$ in $K\{X\}$. Then, set $I(f):=\{g(X)\in K\{X\}: s_{f}^k.g\in \langle f\rangle \}.$ 
\par Let $\theta(\bar x)$ be a quantifier-free $\L_{D}$-formula. Then denote by $\theta^*(\bar y)$ the quantifier-free $\L$-formula obtained by replacing in $\theta(\bar x)$ every term of the form $x_{i}^{(j)}$ by a new variable $y_{ij}$ \cite[Definition 3.16]{GP}.
\enota
\par In this section, $\U$ is a saturated model of $T_{c,D}^*$ containing $K$ and $\A$ as before a saturated model of $DCF_{0}$ containing $\U$. 
(Recall that in a complete theory $T$, if the cardinal $\kappa$ is uncountable and strongly inaccessible, then there is a saturated model of $T$ of cardinality $\kappa$ \cite[Corollary 4.3.14]{Ma}.)
\par To stress that we consider strongly normal extensions within the class of models of $T_{D}$, we will use the term: {\it $T$-strongly normal}.
\dfn Let $K\models T_{D}$ and $L$ be a strongly normal extension of $K$. If $L$ is in addition a model of $T_{D}$, then $L$ is called a $T$-strongly normal extension of $K$.
\edfn
\par Note that by a recent result of Pogudin \cite{Pog}, any finitely generated differentially algebraic field extension $L$ of $K$ such that $C_L\neq L$, is generated by one element; this extends in characteristic $0$ a former result of Kolchin \cite[Chapter 2, Proposition 9]{K73}. Therefore we may assume that any strongly normal extension $L/K$ is of the form $K\langle a\rangle$ for a single element $a$. When giving a criterium below for the existence of $T$-strongly normal extensions, our proof will use axiom scheme (DL) together with this result of Pogudin.

\lem \label{KP} Let $K\models T_{D}$. Let $F/K$ be a strongly normal extension, included in $\A$ and assume $F\neq K$; let 
$a\in \A$ be such that $F=K\langle a\rangle$. Let $\phi(x)$ be a formula isolating $tp^{\A}(a/K)$. 
If there is $\bar b\in \U$ such that $\phi^*(\bar
b)$ holds, then $K$ has a $T$-strongly normal extension in $\U$, isomorphic to $F$ over $K$.
\elem
\pr  By \cite[Lemma 1.4]{M}, $\I^D_{K}(a)$ is of the form $I(f)$ for some $f\in K\{X\}$. Then we may assume that $\phi(x)$ is of the form $f(x)=0\;\&\;s_{f}.q(x)\neq 0$, where $q(X)\in K\{X\}$.
\par Since $f^*(\bar b)=0\;\&\;s_f^*(\bar b)\neq 0$, one can apply the scheme
$(DL)$. Therefore, one can find $u\in \U$ such that $f(u)=0$, with $u^{\n}$ close to $\bar b$. So, we also have that $s_f.q(u)\neq 0$. So $\phi(u)$ holds and by Fact \ref{Tsn}, $K\langle u\rangle$ is a strongly normal extension of $K$ included now in $\U$ and so a model of $T_{D}$.
\qed
\medskip
\par Let us relate the existence result in the above Lemma with the results of M. Kamensky and A. Pillay \cite[Theorem 1.5]{KP}. Let $G$ a connected algebraic group over $C_{K}$ and let $A\in LG(K)$, where $LG$ denotes the Lie algebra of $G$. If $C_{K}$ is existentially closed in $K$, they show that there is a strongly normal extension $L$ of $K$ for the logarithmic differential equation $dlog_{G}(Y)=A$ \cite[Theorem 1.3]{KP}. Furthermore if $C_{K}$ is large and is bounded (namely has for each $n$ only finitely many extensions of degree $n$, the so-called Serre property), then $C_{K}$ is existentially closed in $L$ \cite[Theorem 1.5]{KP}. So, if we apply the results of Kamensky and Pillay when $C_{K}\models T_{c}$ and $C_{K}$ bounded, then automatically $C_{K}$ is large and existentially closed in $K$ and so $K$ has a strongly normal extension $L$ where $C_K$ is existentially closed. That last property implies, since $T$ is universal, that $L\models T_{D}$.

\medskip
\par In Proposition \ref{normalext}, we assume that $dcl^{\U} K\cap L=K$ and $\U$ sufficiently saturated. Note that if $T_{c}=RCF$ and $dcl^{\U}(K)=K$, then $K\models T_{c}$.
\par We will denote by $\U\restriction \L$ the
$\L$-reduct of $\U$. In the following lemma that will be used in Proposition \ref{inter}, we
relate the algebraic closure in models of $T_{c,D}^*$ and of $T_{c}$.  
 \medskip 
\lem \label{acl} Let $L\models T_{D}$ and let $\U$ be a model of $T_{c,D}^*$
extending $L$. Then the algebraic closure $acl^{\U}(L)$ is equal to
$acl^{\U\restriction \L}(L)$.
\elem
\pr Let $a\in acl^{\U}(L)$ and let $\phi(x)$ be an
$\L_{D}(L)$-formula such that $\phi(a)$ holds in $\U$ and which has
only finitely many realizations. Since $T_{c,D}^*$ admits quantifier
elimination, $\phi(x)$ is equivalent to a finite disjunction of
formulas of the form: $$\bigwedge_{i\in I}
p_{i}(x)=0\;\&\;\theta(x),$$ where $\theta^*(\U)$ is an open subset of
some cartesian product of $\U$.
If $I=\emptyset$, then we obtain a contradiction since $\theta(\U)$ is infinite (it is a direct consequence of the scheme (DL) that near every tuple, one can find a tuple of the form $d^{\n}$ \cite[Lemma 3.12]{GP}).
Therefore we may assume that $I\neq \emptyset$; consider $\I^D_{L}(a)$. This is a prime
ideal of the form $\I(f)$, for some $f\in L\{X\}$ \cite[Lemma
1.4]{M}. Note that $f^*(a^{\n})=0\;\&\;s_f^*(a^{\n})\neq 0$. Either the set $S$ of solutions in $\U$ of the formula $f^*(\bar y)=0\;\&\;\theta^*(\bar y)$ is finite and so $a\in acl_{\L}(L)$. Or $S$ is infinite and for any $a\in S$, there exists a neighbourhood $V_{a}$ of $a$, included in $\theta^*(\U)$. Applying the scheme $(DL)$, there exist element $b\in\U\cap V_{a}$ satisfying $f(b)=0$ (and so $\theta(b)$ holds). So we get a contradiction with the finiteness of the number of solutions of $\phi(x)$.  
\qed
\medskip

 \subsection{Picard-Vessiot extensions}
 \par Let $K$ be a differential field which is a model of $T_{D}$ and assume that $C_{K}\models T_{c}$.
\par Let $P(Y):=Y^{(n)}+Y^{(n-1)}.a_{n-1}+\cdots+Y.a_{0}$,  be a differential linear polynomial of order $n$ with coefficients in $K$.
 \dfn \cite[Definition 3.2]{M} Let $L$ be an extension of $K$. Then $L$ is a Picard-Vessiot extension of $K$ corresponding to the linear differential equation $P(y)=0$ if
 \begin{enumerate}
 \item $L$ is generated by the solutions of $P(y)=0$ in $L$,
 \item $C_{L}=C_{K}$,
 \item $P(y)=0$ has $n$ solutions in $L$ which are linearly independent over $C_{K}$.
 \end{enumerate}
 \edfn
 \par Recall that if $L$ be a Picard-Vessiot extension of $K$, then $L$ is a strongly normal \cite[section 9]{M}.
 \par If $C_{K}$ is algebraically closed, there exists a Picard-Vessiot extension of $K$ corresponding to the equation $P(y)=0$, and such extension is unique up to $K$-isomorphism \cite[Theorems 3.4, 3.5]{Magid}. 
 \par In our particular setting, we have the following result.
 The corresponding algebraic equation is $Y_{n}+Y_{n-1}.a_{n-1}+\cdots+Y_{0}.a_{0}=0;$ it has $n$ linearly $C_{K}$-independent solutions in $K$: $u_{i}:=(u_{i0},\cdots, u_{in}),\;i=1,\cdots,n$, equivalently the determinant of the corresponding $n\times n$-matrix formed by the first $n$ coordinates of the $u_{1},\cdots,u_{n}$ is non zero. Let $K\subset \U$ with $\U\models T_{c,D}^*$. By the scheme of axioms $(DL)$, we can find $n$ solutions: $v_{1},\cdots, v_{n}\in \U$ of the differential equation $P(y)=0$ with $(v_{i}, v_{i}^{(1)},\cdots,v_{i}^{(n-1)})$ in any chosen neighbourhood of $(u_{i 0},\cdots, u_{i n-1})$. Moreover by choosing a small enough neighbourhood, we can guarantee that the Wronskian of $v_{1},\cdots, v_{n}$ is non zero, equivalently $v_{1},\cdots, v_{n}$ are linearly independent over $C_{\U}$ \cite[Lemma 4.1]{M}. However in proceeding in this way we did not control the field of constants of $K(v_{1},\cdots,v_{n})$ and so it might be bigger than $C_{K}$. 
\rem A corollary of the proof of  \cite[Theorems 3.4, 3.5]{Magid} is the following. Assume that $C_{K}$ is existentially closed in the class of models of $T$ and let $S$ be the full universal solution algebra for $P$ \cite[Definition 2.12]{Magid}. If there is a maximal differential ideal $I$ of $S$ such that $Frac(S/I)$, the fraction field of $S/I$, is a model of $T_{D }$, then there exists a Picard-Vessiot extension of $K$ corresponding to the equation $P(y)=0$ \cite[Corollary 1.18]{Magid} (which will be $T$-strongly normal).
\erem
\medskip
\lem {\rm  \cite[Theorem 5.5]{Kap}} Let $L$ be a Picard-Vessiot extension of $K$, then $gal(L/K)$ is isomorphic to the $C_{K}$-points of a  linear  algebraic group defined over $C_{K}$. \qed
\elem
\medskip
\ex\label{proper}
When $T$ is the theory of ordered fields, we provide here an example of extensions $K\subset E\subset L$ satisfying the hypotheses of 
	Proposition \ref{interm} and such that $Gal_{\U}(L/K)/Gal_{\U}(L/E)$ is isomorphic to a proper subgroup of $Gal_{\U}(E/K)$.
	
	Let $K:=\R$ endowed with the trivial derivation,
	$E:=\R \langle t\rangle$ where $D(t)=t$,
	$L:=\R \langle u\rangle$ where $u^2=t$ (and so $D(u)=\frac{1}{2}u$),
	$L/K$ and $E/K$ are Picard-Vessiot and so are strongly normal.
	Take $\sigma\in Gal_{\U}(E/K)$ such that $\sigma(t)=-t$. 
	Since $t$ is a square in $L$ and $-t$ is not, $\sigma$ may clearly not be in the image of the map $f$ from the proof of Proposition \ref{interm}.
	
	Here we can easily describe $Gal_{\U}(E/K)$ and $Gal_{\U}(L/K)$. They are both isomorphic to $\mathbb G_m(C_\U)$. 
	
	However  $Gal_{\U}(L/E)$ is isomorphic to the subgroup $\{-1,1\}$ of $\mathbb G_m(C_\U)$. 
	Note that since $C_\U$ is real closed, any element of the group $\mathbb G_m(C_\U)/\{-1,1\}$ is a square, 
	which is no longer true for the group $\mathbb G_m(C_\U)$ and so they cannot be isomorphic.
\eex
\par Note that this example shows that any definable subgroup $G$ of $\mathbb G_m(C_\U)$ cannot be realized as $Gal_{\U}(L/F)$ for $F$ an intermediate extension. Take $G=\mathbb G_m(C_\U)^2.$
\subsection{Weierstrass extensions}\label{W}
Now we will consider the case of a strongly normal extension generated by a solution of a differential equation $W_{k}(Y)=0$ where $W_{k}(Y):= (Y^{(1)})^2-k^2.(4.Y^3+g_{2}.Y-g_{3})\;(\star)$, with $k,\;g_{2},\;g_{3}\in K$ and $27g_{3}^2-g_{2}^3\neq 0$. 
Note that $s_{W_{k}}=\frac{\partial}{\partial Y^{(1)}} W_{k}=2.Y^{(1)}$.
Let $K$ be a differential field satisfying $T_{D}$ and let $\U$ be a model of $T_{c,D}^*$ containing $K$. Then given any solution $(u,u_{1})\in \U^2$ of the Weierstrass equation: $(\dagger)$ $Y^2=4.X^3+g_{2}.X-g_{3}$ with $u_{1}\neq 0$, there exists a solution $a$ of the differential equation $W_{k}(Y)=0$ such that $(a,\frac{a^{(1)}}{k})$ is close to $(u,u_{1})$. 
\par Recall that  an elliptic curve (over $K$) is defined as the projective closure of a nonsingular curve defined by the equation $(\dagger)$. We will denote by $\E$ the corresponding elliptic curve (or $\E(F)$ if we look at its points on an intermediate subfield $F$: $K\subset F\subset \A$) and usually we will use the affine coordinates of the points on $\E$.
We denote by $\oplus$ the group operation on the elliptic curve $\E$ and by $\ominus$ the operation of adding the inverse of an element. It is well-known that $(\E,\oplus)$ is an algebraic group over $K$ (one can express the sum of two elements as a rational function of the coordinates of each of them \cite[Chapter 3, section 2]{Sil}).
\par Let $\si\in Hom_{K}(L,\A)$ and let $a$ be such that $W_{k}(a)=0$. Then $\si(a)$ also satisfies $W_{k}(y)=0$. 
\par Consider the element $(\si(a),\frac{\si(a)^{(1)}}{k}) \ominus (a,\frac{a^{(1)}}{k})$, then one can verify that it belongs to $\E(C_{\A})$ (see \cite[Lemma 2, chapter 3]{K53}, or for instance \cite[section 9]{M}).
So we get the following Lemma.
\lem Let $L:=K\langle a\rangle$, where $a\in \U$ satisfies $W_{k}(a)=0$. Assume that $C_{K}=C_{L}$. Then $L$ is a  strongly normal extension of $K$.\qed
\elem
\par Assume that $L$ is a strongly normal extension of $K$ generated by $a$ with $W_{k}(a)=0$ and set $\mu: gal(L/K)\rightarrow \E(C_{K}): \si\rightarrow (\si(a),\frac{\si(a)^{(1)}}{k}) \ominus (a,\frac{a^{(1)}}{k})$ (respectively $\mu: Gal_{\U}(L/K)\rightarrow \E(C_{\U}): \si\rightarrow (\si(a),\frac{\si(a)^{(1)}}{k}) \ominus (a,\frac{a^{(1)}}{k})$). 
\medskip
\lem \rm{\cite[section 9, Example]{M}} The group $(gal(L/K),\circ)$ \\(respectively $(Gal_{\U}(L/K),\circ)$) is isomorphic to a definable subgroup of  $(\E(C_{K}),\oplus)$ (respectively $(\E(C_{\U}),\oplus)$).
\elem 
\pr It is easy to see that $\mu$ is an injective group morphism \cite[Chapter 3, section 6]{K53}. In order to show that the image of $\mu$ is a definable subset of $\E(C_{K})$ (respectively $\E(C_{\U})$), one notes that $(c_{1},c_{2})\in \E(C_{K})$ belongs to the image of $\mu$ if the first coordinate of $(c_{1},c_{2})\oplus(a,\frac{a^{(1)}}{k})$ satisfies the formula $\chi(x)$ isolating $tp^{\A}(a/K)$. \qed
\medskip
 \par In order to identify $\mu(gal(L/K))$, we now examine in more details which kind of (definable) subgroups may occur in $\E(C_{K})$, when $C_{K}$ is first a real closed field and then a $p$-adically closed field.
\par Let $C_{K}\models RCF$. Then the o-minimal dimension of $\E(C_{K})$ assuming $\E(C_{K})$ infinite, is equal to $1$ and so in case $gal(L/K)$ is infinite, the o-minimal dimension of $\mu(gal(L/K))$ is also equal to $1$. It has a connected component of finite index $\mu(gal(L/K))^{0}$ which is definable \cite[Proposition 2.12]{P88}. So in case $C_{K}=\IR$, we get that $\mu(gal(L/K))^{0}$ is isomorphic to $\IR/\IZ$, as a connected $1$-dimensional compact commutative Lie group over $\IR$ .
\par Recall that over $\IR$, $\E(\IR)$ is either isomorphic to $S^1\times \IZ/2\IZ$ or to $S^1$, depending on the sign of the discriminant (or on whether the polynomial $4.X^3+g_{2}.X-g_{3})$ has three or one real roots) (see \cite{Sil2} or \cite{Poonen}). So the only non-trivial proper subgroup of finite index in $\E(\IR)$ has index $2$; this property (restricted to definable subgroups) can be expressed in a first-order way. So it transfers to any other real closed field, and so $\mu(gal(L/K))^{0}\cong \E(C_{K})^{0}$. Moreover we have that $\mu(gal(L/K))^{0}/\mu(gal(L/K))^{00}\cong \E(C_{K})^{0}/\E(C_{K})^{00}$ where $\mu(gal(L/K))^{00}$ is the smallest type-definable subgroup of bounded index of $\mu(gal(L/K)).$ 
\par Suppose now that $C_{K}$ is a $p$-adically closed field, then one can find a description of definable (and type-definable) subgroups of $\E(C_{K})$ in \cite[Chapter 7]{Sil} and \cite{OP}. Even though in this case the theory does not admit elimination of imaginaries, since we obtained the isomorphism explicitly, we get that $gal(L/K)$ (respectively $Gal_{\U}(L/K)$) is isomorphic to the $C_{K}$-points (respectively to its $C_{\U}$-points) of a group $\L_{rings}(C_{K})$-definable. However, even in the case $C_{K}=\IQ_{p}$, the lattice of subgroups is more complicated. Let us give a quick review of some of the known subgroups. \par First let us recall some notations. Given an elliptic curve $\E$, $\tilde \E$ is its reduction modulo $p$ and $\tilde \E_{ns}$ the subset of the non-singular points of $\tilde \E$. One can also attach to $\E$ a formal group $\hat{\E}$ over $\IZ_{p}$.
Then we have $\E_{1}(\IQ_{p})\subset \E_{0}(\IQ_{p})\subset \E(\IQ_{p})$ and $\E_{1}(\IQ_{p})$ has a subgroup of finite index isomorphic to the additive group of $\IZ_{p}$, where $\E_{0}(\IQ_{p}):=\{P\in \E(\IQ_{p}): \tilde P\in \tilde \E_{ns}(\IF_{p})\}$ and $\E_{1}(\IQ_{p}):=\{P\in \E(\IQ_{p}):\tilde P=\tilde 0\}\cong \hat{\E}(p.\IZ_{p})$ \cite[Proposition 6.3]{Sil}. The index $[\E(\IQ_{p}):\E_{0}(\IQ_{p})]=4$ and $\E_{0}(\IQ_{p})/\E_{1}(\IQ_{p})\cong \tilde \E_{ns}(\IF_{p})$.
\par In case $C_{K}$ is sufficiently saturated, we have $\E(C_{K})$ has an open definable subgroup $H$ such that $H/H^{\circ\circ}$ is isomorphic to the profinite group $(p^r \IZ_{p},+,0)$, where $r>\frac{v(p)}{p-1}$ \cite[Chapter 4, 6.4, b)]{Sil}.
\section{Ordered differential fields}
\par We will use the same notations as in the previous section. We will consider the case where $T$ is the theory of ordered fields and so $T_{c}$ is the theory $RCF$ of real closed fields and $T_{c,D}^*$ the theory $CODF$ of closed ordered differential fields. 
Note that in this case, the language $\La_{D}$ is equal to the language of ordered differential fields, namely $\La_{D}=\{+,-,.,{}^{-1},0,1,<,D\}$. We will use the fact that $CODF$ admits elimination of imaginaries \cite{Point}. 
\par Let $K\subset L\subset \U$ with $K, L\models T_{D}$ and $\U$ a $\vert K\vert^+$-saturated model of $CODF$. Assume that 
$L$ is a strongly normal extension of $K$. Recall that  
both $gal(L/K)$ (respectively $Gal_{\U}(L/K)$) are isomorphic to the $C_{K}$-points (respectively the $C_{\U}$-points) of an algebraic group  $\L_{rings}(C_{K})$-definable. Since $\U\models CODF$, $C_{\U}$ is a model of $RCF$. Given a subfield $F$ of $\U$, we will denote by $F^{rc}$ its real-closure in $\U$.
\par By former results of A. Pillay \cite{P88}, any definable group $G$ in a real closed field can be endowed with a definable topology in such a way it becomes a topological group. Moreover $G$ has a connected component $G^0$ of finite index (which is also definable) \cite[Proposition 2.12]{P88}. If $G$ lives in a sufficiently saturated model, then it has a smallest type-definable subgroup of bounded index denoted by $G^{00}$; this result holds in fact in any $NIP$ theory \cite{Sh}. One can endow the quotient $G/G^{00}$ with the logic topology and if $G$ is in addition definably connected and definably compact, then $G/G^{00}$ is a compact connected Lie group \cite[Theorem 1.1]{BOPP}. Moreover the $o$-minimal dimension of $G$ is equal to the dimension $G/G^{00}$ as a Lie group \cite[Theorem 8.1]{HPP}.
\nota Given $F_0\subset F_1$ be two subfields of $\U$. We denote by $aut(F_1/F_0)$ the group of $\L_{D}$-automorphisms of $F_1$ fixing $F_0$ and by $Aut_{T}(L/K)$ the group \\$aut(\langle L, C_{\U}\rangle/\langle K, C_{\U}\rangle)$.
\enota
\par Note that any $\si\in Aut_T(L/K)$ has a unique extension to an element of \\
$aut(\langle L,C_{\U}\rangle^{rc}/\langle K,C_{\U}\rangle)$.

In the next proposition
we will use the following easy consequence of Lemma \ref{acl} that
$\langle L,C_{\U}\rangle^{rc}=dcl^{\U}(\langle L,C_{\U}\rangle)$ (one uses in addition that in an ordered structure,
the algebraic closure is equal to the definable closure).
\par Let $L$ be a strongly normal extension of $K$ and $\bar d\in \langle L,C_{\U}\rangle^{rc}$. 
Then $L\langle\bar d\rangle$ is again a strongly normal extension of $K\langle\bar d\rangle$.
(Indeed, $C_{L\langle \bar d\rangle}=C_{L}=C_{K}$ since $C_{K}=C_{L}$ is algebraically closed in $C_{\U}$ and if $\si \in Hom_{K\langle\bar d\rangle}(L\langle\bar d\rangle,\A)$, then $\si(L\langle\bar d\rangle)\subset \langle L\langle\bar d\rangle,C_{\A}\rangle$.) 
If $\si\in Hom_{K\langle\bar d\rangle}(L\langle\bar d\rangle,\U)$, then $\si(L)\subset \langle L,C_{\U}\rangle$ and $\si(\bar d)=\bar d$. So, $\si(L\langle\bar d\rangle)\subset \langle L\langle\bar d\rangle,C_{\U}\rangle$.
\par Since $\bar d\in \langle L,C_{\U}\rangle^{rc}=dcl^{\U}(\langle L,C_{\U}\rangle)$, any $\tau\in Aut_{T}(L/K)$ extends uniquely to an element of $Aut_{T}(L\langle\bar d\rangle/K\langle\bar d\rangle)$. 
 \par Note that $\mu(Aut_{T}(L/K))$ (see Notation \ref{mu}) is a subgroup of $H(C_\U)$ which is an infinitely $\L_D(L)$-definable group in $\U$ (one expresses that the image of any positive element of $\langle L,C_{\U}\rangle$ remains positive). A natural question which we have not answered is whether $\mu(Aut_{T}(L/K))$ is already $\L_D(L)$-definable. We consider intersections of $\mu(Aut_{T}(L/K))$ with groups of the form $\mu(G_0)$, with $G_0$ a subgroup of $Gal_{\U}(L/K)$ such that $\mu(G_0)$ is $\L_D(\langle L,C_{\U}\rangle)$-definable subgroup in $\U$; we will call such intersections {\it relatively $\L_D(\langle L,C_{\U}\rangle)$-definable subgroups of $\mu(Aut_T(L/K))$}. 
\prop\label{inter} Let $\U$ be a $\vert K\vert^+$-saturated model of $CODF$ and let $L/K$ be a strongly normal extension with $C_L\subset_{ac} \U$.
\par Let $G_0$ be a subgroup of $Gal_{\U}(L/K)$ and assume that $\mu(G_0)$ is $\L_D(\langle L,C_{\U}\rangle)$-definable in $\U$. Then there is a tuple of elements $\bar e\in \langle L,C_{\U}\rangle^{rc}$ such that we have: $$G_{0}\cap Aut_{T}(L/K)\cong Aut_T(L\langle \bar e\rangle)/K\langle \bar e\rangle).$$
\eprop
\pr  Let $L:=K\langle a\rangle$ for some $a\in L$. As in the proof of Theorem \ref{Kolchin}, for $\si\in Gal_{\U}(L/K)$, we write $\si(a)$ as $\tilde f(a,\bar d, \bar k)$
with $\bar k\in K$, $\bar d\in C_{\U}$.
Let $G_{0}$ be a subgroup of $Gal_{\U}(L/K)$ such that $\mu(G_0)$ is definable (in $\U$) by some $\L_D(\langle L,C_{\U}\rangle)$-formula $\psi(\bar x)$, namely $\si\in G_0$ if and only if $\psi(\bar d)$ holds. 

\par Consider $W:=\{\tau(a):\;\tau\in G_{0}\}$, it is a $\L_{D}(\langle L,C_{\U}\rangle)$-definable subset of $\langle L,C_{\U}\rangle$. 
So $W$ is an $\L_{D}(\langle L,C_{\U}\rangle)$-definable subset of $\U$.
\par Since $CODF$ has elimination of imaginaries \cite[Theorem 2.6]{Point}, there is a canonical parameter $\bar e$ for $W$ belonging to $dcl^{\U}(\langle L,C_{\U}\rangle)=\langle L,C_{\U}\rangle^{rc}$ (Lemma \ref{acl}).
\par Set $F:=K\langle \bar e\rangle$. Let us show that $G_{0}\cap Aut_{T}(L/ K)\cong Aut_{T}(L\langle\bar e\rangle/F)$.
\par Let $\tau\in Aut_{T}(L\langle\bar e\rangle/F)$ and consider $\tau\restriction L\langle\bar e\rangle$. 
Since $\U$ is $\vert K\vert^+$-saturated, it extends to an automorphism $\tilde \tau$ of $\U$, which leaves $F$ fixed and so it leaves $W$ invariant.  Since $a\in W$, $\tau(a)\in W$, namely there is $\tau_{0}\in G_{0}$ and $\tau(a)=\tau_{0}(a)$.
\par Since $L\langle \bar e\rangle$ is a strongly normal extension of $F$ and $\tau(L\langle \bar e\rangle)\subset \U$, $\langle L,C_{\U}\rangle=\langle \tilde \tau(L),C_{\U}\rangle=\langle \tau(L),C_{\U}\rangle$ (see Remark \ref{star}). 
In particular, $\tilde\tau(a)=\tau(a)\in \langle L,C_{\U}\rangle$. 

\par Conversely, let $\tau_0\in G_{0}\cap Aut_{T}(L/K)$ and $\tilde \U$ be a $\vert \langle L,C_{\U}\rangle\vert^+$-saturated extension of $\U$. Extend $\tau_0$ to an $\La_{D}$-automorphism $\si$ of $\tilde\U$ which is the identity on $\langle K, C_{\U}\rangle$. Since $\si(a)=\tau_0(a)$, $\si$ leaves $W$ invariant and thus $\bar e$ fixed.
 Therefore the restriction of $\si$ to $\langle L\langle\bar e\rangle, C_{\U}\rangle$ belongs to $Aut_{T}(L\langle\bar e\rangle/ F)$. 
 \qed
\medskip
\par Note that we have obtained a Galois correspondence between relatively $\L_D(\langle L,C_{\U}\rangle)$-definable subgroups of $\mu(Aut_{T}(L/K))$  and intermediate finitely generated differential field extensions of $K$ in $L^{rc}$.
Indeed, by Proposition \ref{inter}, given $G_0$ a subgroup of $Gal_{\U}(L/K)$ such that $\mu(G_0)$ is $\L_D(\langle L,C_{\U}\rangle)$-definable, we can find a tuple $\bar e\in \langle L,C_{\U}\rangle^{rc}$ such that $G_{0}\cap Aut_{T}(L/K)\cong Aut_{T}(\langle L,C_{\U}\rangle^{rc}/K\langle \bar e\rangle)$. Conversely given any tuple $\bar e\in \langle L,C_{\U}\rangle^{rc}$, we can find $G_0\subset Gal_{\U}(L/K)$ such that $Aut_{T}(\langle L,C_{\U}\rangle^{rc}/K\langle \bar e\rangle)$ is 
of the form $G_{0}\cap Aut_{T}(\langle L,C_{\U}\rangle^{rc}/K)$. Indeed since $L$ is a strongly normal extension of $K$, given any $\tau \in Aut_{T}(\langle L,C_{\U}\rangle^{rc}/K)$, $\tau(\langle L, C_{\U}\rangle)\subset \langle L, C_{\U}\rangle$. Given each component $e_i$ of $\bar e$ and its minimal polynomial $p_i(x)$ with coefficients in $\langle L, C_{\U}\rangle$, we define $G_0$ as the subset of $Gal_{\U}(L/K)$ leaving fixed the coefficients of each $p_i(x)$. To see that $\mu(G_0)$ is  $\L_D(\langle L,C_{\U}\rangle)$-definable, we proceed as in Corollary \ref{GalE}. We use the facts that each coefficients of the $p_i(x)$ can be expressed as $h(a)$ for some  $\L_D(\langle K,C_{\U}\rangle)$-definable function and that for $\tau\in Gal_{\U}(L/K)$, $\tau(a)=\tilde f(a,c_{a,\tau(a)},\bar k)$.


\section{ On the non finitely generated case}\label{nfg}
J. Kovacic extended the notion of strongly normal extensions to the infinitely generated case \cite[Chapter 2, section 1]{Kov1}. Recall that $L$ is a strongly normal extension of $K$ if $L$ is a union of finitely generated strongly normal extensions of $K$. Equivalently, $C_{L}=C_{K}$ and $L$ is generated by finitely generated strongly normal extensions of $K$ \cite[Chapter 2, Section 1, Definition]{Kov1}. 
\par In this section we will place ourselves in the setting of section 3. We will fix a differential field $K$, we work inside an extension $\U$, and embed $\U$ in a saturated model $\A$ of $DCF_0$.
\par In \cite[Chapter 2, Theorem 1]{Kov1}, J. Kovacic proved that if $L$ is a strongly normal extension of $K$ with $C_{K}$ algebraically closed, then on one hand, the set of finitely generated strongly normal extensions of $K$ in $L$ forms an injective system and on the other hand, $Gal(L/K)$ is the inverse limit of the differential Galois groups of finitely generated extensions. In particular, it is isomorphic to an inverse limit of algebraic groups.
\par In the following proposition, we will establish an analog for the differential Galois groups $gal(L/K)$ and $Gal_{\U}(L/K)$ (Notation \ref{auto-nota}) (without the assumption that $C_K$ is algebraically closed). 

\medskip
\prop Let $L/K$ be a strongly normal extension with $L\subset \U$.
Then $Gal_{\U}(L/K)$) is isomorphic to a subgroup of an inverse limit of groups of the form $H_i(C_{\U})$ where $H_i$ is an algebraic group defined over $C_{K}$, $i\in I$. Moreover, $gal(L/K)$ is a subgroup of $Gal_{\U}(L/K)$ and is isomorphic to a subgroup of an inverse limit of groups of the form $H_i(C_K)$, $i\in I$.
\eprop 
\pr Let $L=\bigcup_{F\in \F} F$, where $\F$ denotes the set of all finitely generated strongly normal extensions of $K$ in $L$. Then $Gal_{\U}(L/F)$ is a normal subgroup of $Gal_{\U}(L/K)$ (see Corollary \ref{GalE}, note that there we only use that $F$ was a (finitely generated) strongly normal extension of $K$). 
Moreover we have an embedding of $Gal_{\U}(L/K)/Gal_{\U}(L/F)$ into $Gal_{\U}(F/K)$, sending $\si$ to $\si\restriction F$ (Corollary \ref{interm}).

\par  For each $K\subset F_{i}\subset F_{j}\subset L$ with $F_{i},\,F_{j}\in \F$, define the maps $$f_{F_{i}}^{F_{j}}:Gal_{\U}(L/K)/Gal_{\U}(L/F_{j})\rightarrow Gal_{\U}(L/K)/Gal_{\U}(L/F_{i}),$$ sending $\si.Gal_{\U}(L/F_{j})$ to $\si.Gal_{\U}(L/F_{i})$ and $$r_{F_{i}}^{F_{j}}:Gal_{\U}(F_{j}/L)\rightarrow Gal_{\U}(F_{i}/K)$$ sending $\si$ to $\si\restriction F_{i}$ (with kernel $Gal_{\U}(F_{j}/F_{i})$).
The following diagram commutes:
\medskip

\begin{center}
$\begin{CD}
Gal_{\U}(L/K)/Gal_{\U}(L/F_{j}) @>>> Gal_{\U}(F_{j}/L)\\
@Vf_{F_{i}}^{F_{j}}VV  @Vr_{F_{i}}^{F_{j}}VV\\
Gal_{\U}(L/K)/Gal_{\U}(L/F_{i}) @>>> Gal_{\U}(F_{i}/K)
\end{CD}$
\end{center}
\medskip
\par The family of differential Galois groups $Gal_{\U}(F/K)$ forms a projective system, as well as the family of quotients $Gal_{\U}(L/K)/Gal_{\U}(L/F)$, $F\in \F$. \\Since $\bigcap_{F\in \F} Gal_{\U}(F/K)=\{1\}$, we get an embedding of $Gal_{\U}(L/K)$ into $\varprojlim_{F\in \F} Gal_{\U}(F/K)$.
\par Finally, for $\si\in Gal_{\U}(L/K)$, define $\mu(\si):=(\mu(\si\restriction F))$ (Notation \ref{mu}).\\
We get $\mu(Gal_{\U}(L/K))\subset \varprojlim_{F\in \F}\mu(Gal_{\U}(F/K))$.
\par The second part is clear. \qed
\medskip
\par Let us give an example of a strongly normal extension which is infinitely generated. 
\ex Let $K$ be a formally real differential field with a real closed subfield of constants.  Let $K^*$ be the intersection of all real closed subfields of a given algebraic closure of $K$. Note that $K^*$ is a pythagorean field (namely a sum of two squares is a square). By a result of M. Griffin \cite{C}, $K^*$ is the maximal Galois extension of $K$ to which all the orderings of $K$ extend \cite[Theorem 1.1]{C}. 
Since any finite Galois extension $L$ of $K$ such that $C_K$ is algebraically closed in $L$ is a Picard-Vessiot extension of $K$ \cite[Propositions 3.9, 3.19, 3.20]{Magid} (the hypothesis there is that the field of constants is algebraically closed, however one can check that it can be weakened to {\it the field of constants is relatively algebraically closed}), we obtain that 
$K^*$ is a union of Picard-Vessiot extensions of $K$. Moreover $K^*$ is infinitely generated over $K$ whenever $K$ is not pythagorean \cite[chapter 9]{R}.
\par Let $F$ be a real closed field and take for $K$ either the formally real differential field $F(t)$ or the Laurent series field $F((t))$, where $D$ is trivial on $F$ and $D(t)=1$. Then neither fields $K$ are pythagorean and so the corresponding $K^*$ are concrete examples of infinitely generated Picard-Vessiot extensions.
\eex
\par By Zorn's lemma, one can always find a maximal strongly normal extension of $K$ inside $\U$. Adapting the classical proofs \cite{Magid}, we will show that one can find a maximal strongly normal extension of $K$ inside $\U$ which has proper strongly normal extension (inside $\U$). 
\par Note that under the assumption that $C_K$ is algebraically closed, taking such extension does not depend on the universal model of $DCF_{0}$ one chooses. Moreover,  one can show that two maximal strongly normal extensions of $K$ in $\A$ are isomorphic over $K$ \cite[Chapter 2, section 4, Corollary]{Kov1}.

\rem Let $K\models T_{c}$ with a trivial derivation and suppose that $K$ is $\aleph_{0}$-saturated (we have to be able to find $n\geq 3$ elements $\alpha_{1},\cdots,\alpha_{n}$ of $K$ linearly independent over $\IQ$). Let $\U$ be a saturated extension of $K$ satisfying $T_{c,D}^*$. Inside $\U$, choose $n$ elements $a_{1},\cdots, a_{n}$ which are algebraically independent over $K$. Consider $L_{2}$ the differential subfield generated by $K$ and $a_{1},\cdots, a_{n}$. Consider the subfield $L_{1}$ of $L_2$ containing $K$ and generated by $b_{1},\cdots, b_{n}$ the coefficients of the polynomial $\prod_{i=1}^n (X-a_{i})$ (namely the $b_{i}$'s are the images (up to a sign) of the elementary symmetric functions on $a_{1},\cdots, a_{n}$ and so also algebraically independent over $K$). 
Now define a derivation on $L_{2}$ by first imposing that $D(b_{i})=\alpha_{i}.b_{i}$ and then $D(a_{i})$ is given by using the fact that $a_{i}$ is algebraic over $L_{1}$.
\par The field $L_{1}$ is a Picard-Vessiot extension of $K$
\cite[Proposition 1.20, Example 3.23]{Magid} corresponding to the linear differential
equation $Y^{(n)}+Y^{(n-1)}.\beta_{n-1}+\cdots+Y.\beta_{1}=0$,
where $\prod_{i=1}^n
(Y-\alpha_{i})=Y^n+Y^{n-1}.\beta_{n-1}+\cdots+Y.\beta_{1}$. The extension
$L_2/L_1$ is a Picard-Vessiot extension. Indeed, we use \cite[Proposition 3.20]{Magid} (with the same proviso as above), 
$L_2$ is a finite
Galois extension of $L_{1}$ and $C_{L_{1}}=C_{K}=K$ is algebraically
closed in $\U$. 
However $L_{2}$ is not contained in a Picard-Vessiot extension of $K$. Suppose otherwise,
then $\langle L_{2},\bar C_{K}\rangle$ would also be contained in a Picard-Vessiot extension of $\langle K, \bar C_{K}\rangle$, which contradicts \cite[Example 3.33]{Magid}.
\par Now let $F$ be a maximal strongly normal extension of $K$ inside $\U$ containing $L_{1}$. Let us show that the subfield $F_{1}:=\langle F, L_{2}\rangle$ of $\U$ is a proper strongly normal extension of $F$. 
\par Indeed, $F_{1}/F$ is a Picard-Vessiot extension since it is a finite
Galois extension and $C_{F_1}=C_K$. Suppose that $F_1=F$, then
$\langle L_{2},\bar C_{K}\rangle$ is
contained in a strongly normal extension of $\langle K, \bar
C_{K}\rangle$. It implies by \cite[Proposition 11.4]{Kov3} that $\langle L_{2},\bar C_{K}\rangle$ is 
contained in a Picard-Vessiot extension of $\langle K,\bar
C_{K}\rangle$, a contradiction.
\erem
\medskip
\par
\noindent {\bf Acknowledgments:}
Both authors thank the MSRI for their hospitality during the spring of 2014 and the FRS-FNRS for their support. 
 They thank Omar Leon Sanchez and an anonymous referee for their kind and useful remarks which led to a significant improvement of the paper. They thank Anand Pillay for his insightful remarks, in particular for insisting on the fact that a previous hypothesis on the field of constants of the intermediate field $\U$ was not necessary, and on Lemma 5.4 (and the comment afterwards). We are also grateful to Michael Singer for his advice.

\end{document}